\documentclass[letterpaper,12pt]{article}

\usepackage[intlimits]{amsmath}
\usepackage{amssymb}
\usepackage{amsthm}
\usepackage{mathrsfs}
\usepackage{cite}
\usepackage{color}
\usepackage{bm}

\newcommand{\pr}{\partial}
\newcommand{\dom}{\Omega}

\newcommand{\spn}{{\mathrm{span}\,}}

\newcommand{\R}{\mathbb{R}}

\newcommand{\Ord}{\mathscr{O}}

\newcommand{\dd}{\,\text{d}}

\newtheorem{thm}{Theorem}

\title{Prescribed nonlinearity helps in an anisotropic Calder\'on-type problem}
\author{C\u{a}t\u{a}lin I. C\^{a}rstea\thanks{Department of Applied Mathematics, National Yang Ming Chiao Tung University, Hsinchu 300, Taiwan, R.O.C.; \mbox{email: catalin.carstea@gmail.com}}}

\date{}

\begin{document}
\maketitle

\begin{abstract}
In this paper I consider the inverse boundary value problem for a quasilinear, anisotropic, elliptic equation of the form $\nabla\cdot(\gamma\nabla u+|\nabla u|^{p-2}\nabla u)=0$, where $\gamma$ is a smooth, matrix valued, function with a uniform lower bound. I show that boundary Dirichlet and Neumann data for this equation, in the form of a Dirichlet-to-Neumann map, determine the coefficient matrix uniquely, in dimension 3 and higher. This stands in contrast to the classical linear anisotropic Calder\'on problem where there is a known obstruction to uniqueness due to the invariance of the boundary data under transformations of the equation via any boundary fixing diffeomorphism. 
\end{abstract}

\section{Introduction}\label{introduction}

Let $\dom\subset\R^n$, $n\geq3$ be a closed, bounded, connected, smooth domain. Consider the boundary value problem 
\begin{equation}\label{eq-general}
\left\{\begin{array}{l}\nabla\cdot (a(x,u, \nabla u)\nabla u)=0,\\[5pt] u|_{\pr\dom}=f,\end{array}\right.
\end{equation}
with a suitably chosen matrix valued function $a(x,p,\xi)$, such that the above equation becomes elliptic. Then one can introduce the Dirichlet-to-Neumann map
\begin{equation}
\Lambda_a(f)=\nu\cdot a(x,u,\nabla u)\nabla u|_{\pr\dom}.
\end{equation}
In \cite{Ca}, in the particular case $a_{jk}(x,p,\xi)=\gamma(x)\delta_{jk}$, Calder\'on asked the question: can $a$ be recovered from $\Lambda_a$? In this linear isotropic case the answer is yes, as first shown in \cite{SyUh} (for $n\geq3$) and \cite{Na} (for $n=2$). In the isotropic quasilinear case $a_{jk}(x,p,\xi)=\gamma(x,p,\xi)\delta_{jk}$, it was shown in \cite{CarFeKiKrUh}, for $n\geq3$ and a large class of functions $\gamma(x,p,\xi)$ that the answer is also yes.

In the generally anisotropic case where $a$ is not assumed to be a diagonal matrix, it has been observed (first by Luc Tartar in the linear case; account given in \cite{KoVo}) that uniqueness cannot hold. If $\Phi: \dom\to \dom$ is a smooth diffeomorphism that acts like the identity on $\pr\dom$, then the coefficient matrix
\begin{equation}\label{eq-tartar}
\tilde a (x,p,\xi)=\frac{1}{D\Phi}(D\Phi)^Ta(\cdot,p,D\Phi\xi)D\Phi\circ\Phi^{-1}(x)
\end{equation}
produces the same Dirichlet-to-Neumann map as $a$. To date, even in the linear case (except in dimension $n=2$), it is not known if this is the only obstruction to uniqueness. If $a_{jk}(x,p,\xi)=\gamma_{jk}(x)$, some results are available, e.g. in \cite{LeeUh} for real analytic coefficient matrices $\gamma$,  in \cite{DoKeSaUh}, \cite{DoKuLaSa} for coefficient matrices that satisfy certain structure constraints, and in \cite{CarFeOk} for an assortment of related questions.

Recently there has been a considerable increase in the study of inverse problems for non-linear equations. Examples of works in this direction include   \cite{FeOk}, \cite{Is1}, \cite{IsNa}, \cite{IsSy},  \cite{KrUh}, \cite{KrUh2}, \cite{LaLiiLinSa1}, \cite{LaLiiLinSa2}, \cite{Sun2} on semilinear equations, and \cite{CarFe1}, \cite{CarFe2}, \cite{CarFeKiKrUh}, \cite{CarGhNa}, \cite{CarGhUh}, \cite{CarKa}, \cite{CarNaVa}, \cite{Car}, \cite{EgPiSc}, \cite{HeSun}, \cite{Is2}, \cite{KaNa}, \cite{MuUh}, \cite{Sh}, \cite{Sun1}, \cite{Sun3}, \cite{SunUh} on quasilinear equations. I also want to mention separately \cite{BrKaSa}, \cite{Br}, \cite{BrHaKaSa},  \cite{BrIlKa}, \cite{CarFe3}, \cite{GuKaSa}, \cite{KaWa}, \cite{SaZh}, which deal with equations involving the $p$-Laplacian and related operatos. Most of these works employ the  so called \emph{second/higher linearization} method, which first appeared in \cite{Is1}. In this method one uses Dirichlet data that depends on a small (or large) parameter $\epsilon$, typically of the form $\epsilon\phi$, then considers the asymptotic expansion in terms of the parameter $\epsilon$ of $\Lambda(\epsilon\phi)$, the Dirichlet-to-Neumann map evaluated on this data, in order to obtain information about the coefficients of the equation. Many times this is presented as differentiating the equation and $\Lambda(\epsilon\phi)$ a number of times with respect to the small parameter, then setting the parameter to zero.

In this paper I will consider a particular case of \eqref{eq-general}. Let $p\in(1,2)\cup(2,\infty)$ and let $\gamma_{jk}\in C^\infty(\dom)$ be the coefficients of an $n\times n$ matrix $\gamma$, which further satisfies the ellipticity condition
\begin{equation}
\lambda|\xi|^2\leq \xi\cdot\gamma(x)\xi\leq \lambda^{-1}|\xi|^2, \quad\forall x\in\dom,\xi\in\R^n.
\end{equation}
The boundary value problem that I would like to consider is
\begin{equation}\label{eq}
\left\{\begin{array}{l}\nabla\cdot (\gamma(x)\nabla u)+\nabla\cdot(|\nabla u|^{p-2}\nabla u)=0,\\[5pt] u|_{\pr\dom}=f.\end{array}\right.
\end{equation}
The forward problem is well posed.
\begin{thm}[e.g. \mbox{\cite[Theorem 1]{Lie}}]\label{thm-existence}
Let $f\in C^{1,\alpha}(\pr\dom)$, for some $\alpha \in (0,1]$. There exist $\beta\in (0,1)$, and $C(\|f\|_{C^{1,\alpha}(\pr\dom)})>0$ nondecreasing, such that equation \eqref{eq} has a unique weak solution $u\in C^{1,\beta}(\dom)$ and
\begin{equation}
\|u\|_{C^{1,\beta}(\dom)}\leq C(\|f\|_{C^{1,\alpha}(\pr\dom)}).
\end{equation}
\end{thm}
\noindent Therefore, the Dirichlet-to-Neumann map associated to \eqref{eq},
\begin{equation}
\Lambda_\gamma(f)=\left.\left(\nu\cdot\gamma\nabla u+|\nabla u|^{p-2}\pr_\nu u  \right)\right|_{\pr\dom},
\end{equation}
can be defined. The main result of this paper is
\begin{thm}\label{thm-main}
 The coefficient matrix $\gamma$ can be explicitly computed from $\Lambda_\gamma$.
\end{thm}
Note that there is no diffeomorphism invariance here. The intuition for this  is that any boundary fixing diffeomorphism that transforms the $\gamma$-equation to a different $\tilde\gamma$-equation would have to separately transform the linear part and the non-linear part, according to \eqref{eq-tartar}. Since the nonlinear part is prescribed apriori, this should require the diffeomorphism to be the identity. One could say that here the nonlinearity helps. This result supports the conjecture that diffeomorphism invariance is the only obstruction to uniqueness in the anisotropic Calder\'on inverse boundary value problem.

Another result in which a prescribed nonlinearity is known to help in the inverse problem for an anisotropic elliptic equation is \cite{CarFeOk}[Theorem 6]. There it is shown that for semilinear equations of form 
\begin{equation}
-\triangle_g u+uG(x,u)=0,
\end{equation}
where $\triangle_g$ is the Laplace-Beltrami operator for a Riemann manifold with boundary $M,g)$, and with the nonlinearity $G$ satisfying certain conditions, the identity the Dirichlet-to-Neumann maps corresponding to metrics $g_1$ and $g_2$ implies the existence of a boundary fixing diffeomorphism $\Phi$ such that $g_1=\Phi^*g_2$, which is the expected obstruction for equations in Laplace-Beltrami form.

Inverse problems for an equation like \eqref{eq} have been studied before numerically in \cite{HaHyMu}, and using a version of the linearization method in \cite{CarKa}. In both works the coefficients being recovered are scalar. There is also a connection to the papers \cite{CarFe1}, \cite{CarFe2}. In these works an inverse problem for a quasilinear anisotropic equation of the form 
\begin{equation}
\triangle u +\nabla \cdot(A:\nabla u\otimes\nabla u)=0,
\end{equation}
 where $A$ is an arbitrary 3-tensor, is considered.  As in the case of the present paper, any boundary fixing diffeomorphism should transform the linear and the quadratic terms separately. Since the linear term is prescribed, the expectation is that the Dirichlet-to-Neumann map determines the tensor $A$ uniquely, which is indeed the result obtained there.

The idea for the proof of Theorem \ref{thm-main} resembles the linearization method, but with a difference. For the sake of simplicity, suppose $p\in(2,\infty)$.  In this case the typical linearization approach would be to take boundary data of the form $f=\epsilon \phi$, with $\epsilon$ a small parameter. Since the corresponding solution $u_\epsilon$ is of size $\epsilon$, the linear term of the equation will also be of size $\epsilon$, while the nonlinear term will be of size $\epsilon^{p-1}$, which makes it negligible in the first approximation. One can prove that, in suitable norms,
\begin{equation}
u_\epsilon=\epsilon u_0+o(\epsilon),\quad \text{ where } \nabla\cdot(\gamma\nabla u_0)=0,\quad u_0|_{\pr\dom}=\phi.
\end{equation}
It would then follow that the nonlinear Dirichlet-to-Neumann map $\Lambda_\gamma$ determines the linear Dirichlet-to-Neumann map for the equation $\nabla\cdot(\gamma\nabla u_0)=0$. If $\gamma$ were scalar, uniqueness for it would follow. This is the approach taken in \cite{CarKa}. (The $p\in(1,2)$ is analogous, with $\epsilon$ replaced with $\epsilon^{-1}$.)

In this paper I take the alternate approach of using boundary data of the form $f=\epsilon^{-1}\phi$ in \eqref{eq}, with $\epsilon$ a small parameter. The expectation then is that the nonlinear term dominates as $\epsilon\to0$, and indeed it is shown in section \ref{large-asymptotics} that the corresponding solution $u_\epsilon$ is of the form
\begin{equation}
u_\epsilon=\epsilon^{-1} v+\epsilon^{p-3} R+o(\epsilon^{p-3}),
\end{equation}
where
\begin{equation}
\left\{\begin{array}{l}\nabla\cdot(|\nabla v|^{p-2}\nabla v)=0,\\ v|_{\pr\dom}=\phi.\end{array}\right.
\end{equation}
The term $R$ satisfies the equation
\begin{equation}
\left\{\begin{array}{l}\nabla\cdot(A(v)\nabla R)=-\nabla\cdot(\gamma\nabla v),\\[5pt] R|_{\pr\dom}=0,\end{array}\right.
\end{equation}
where $A(v)$ is the matrix with coefficients
\begin{equation}\label{eq-def-A}
A(v)_{jk}=|\nabla v|^{p-2}\left(\delta_{jk}+(p-2)\frac{\pr_j v\pr_k v}{|\nabla v|^2}\right).
\end{equation}
As it will be clear from the proof, and is also suggested by the definition of $A(v)$, it is necessary to choose $v$ so that $\nabla v$ never vanishes. Fortunately, there are many choices, e.g. any linear function of $x$ would work. Expanding $\Lambda_\gamma(\epsilon\phi)$ in powers of $\epsilon$, it will follow that $\Lambda_\gamma$ determines
\begin{equation}
I(v,W)=\int_\dom\nabla v\cdot\gamma\nabla W-\int_\dom R\nabla\cdot\left(A(v)\nabla W\right), \quad\forall W\in C^\infty(\dom).
\end{equation}
This turns out to be sufficient to determine $\gamma$, with appropriate choices for $v$ and $W$. These choices involve the linearization of the $p$-Laplace equation itself, using boundary data of the form $\phi=\phi_0+\tau\phi_1$ depeding on a parameter $\tau$. In section \ref{proof-main-thm}, by differentiating with respect to the parameter $\tau$, I derive an integral identity involving $\gamma$ and two arbitrary solutions of the equation $\nabla\cdot(A(v_0)\nabla V)=0$, with $v_0$ a linear function of $x$. Choosing solutions inspired, in two different ways, by the ones proposed by Calder\'on in \cite{Ca}, I can then show that the Fourier transform $\hat\gamma$ can be computed from $\Lambda_\gamma$.

In section \ref{p-Laplace-linearization} I detail the linearization procedure for the $p$-Laplace equation, for the whole rage of $p\in(1,2)\cup(2,\infty)$. In section \ref{small-asymptotics} I present the small data asymptotics for equation \eqref{eq}, when $p\in(1,2)$ and in section \ref{large-asymptotics} I present the large data asymptotics when $p\in(2,\infty)$. Finally, section \ref{proof-main-thm} gives the proof of Theorem \ref{thm-main}.

\section{Linearization of the $p$-Laplace equation}\label{p-Laplace-linearization}

This section is concerned with a few observations about the boundary value problem for the $p$-Laplace equation
\begin{equation}\label{eq-p-Laplace}
\left\{\begin{array}{l}\nabla\cdot(|\nabla v|^{p-2}\nabla v)=0,\\[5pt] v|_{\pr\dom}=f.\end{array}\right.
\end{equation}
By e.g. \cite[Theorem 1]{Lie}, if $f\in C^{1,\alpha}(\pr\dom)$, then there exists a $\beta\in (0,1)$ so that a unique weak solution $v\in C^{1,\beta}(\dom)$ exists, and $||v||_{C^{1,\beta}(\dom)}\leq C(||f||_{C^{1,\alpha}(\pr\dom)})$. 

Equation \eqref{eq-p-Laplace} has solutions that have no critical points in $\dom$. For example, any coordinate function $x_j$ is such a solution. Suppose $\phi_0\in C^\infty(\pr\dom)$ is the boundary value of one such solution $v_0\in C^\infty(\dom)$. If $\phi_1\in C^\infty(\pr\dom)$ is arbitrary and $-1<\tau<1$ is a small parameter, let $v_\tau$ be the solution to 
\begin{equation}\label{eq-tau}
\left\{\begin{array}{l}\nabla\cdot(|\nabla v_\tau|^{p-2}\nabla v_\tau)=0,\\[5pt] v_\tau|_{\pr\dom}=\phi_0+\tau\phi_1.\end{array}\right.
\end{equation}
Since their boundary data is bounded, the solutions $v_\tau$ are uniformly bounded in $C^{1,\beta}(\dom)$. By the theorem of Arzel\'a-Ascoli, there must exist $\tilde v_0\in C^{1,\beta}(\dom)$ such that (after passing to a subsequence) $v_\tau\to v$ in $C^1(\dom)$. Taking the limit of \eqref{eq-tau}, in the sense of $\mathscr{D}'(\dom)$, it follows that $v$ solves \eqref{eq-p-Laplace}, with Dirichlet data $v|_{\pr\dom}=\phi_0$. This implies that $v=v_0$, and therefore $v_\tau-v_0\to 0$ in $C^1(\dom)$. Incidentally, this shows that $v_\tau$, for $\tau$ small enough, also has no critical points. 

Let $V_\tau$ be defined by the Ansatz
\begin{equation}
v_\tau=v_0+\tau V_\tau.
\end{equation}
For $\xi\in\R^n\setminus\{0\}$ let
\begin{equation}
J_j(\xi)=|\xi|^{p-2}\xi_j,\quad j=1,\ldots,n.
\end{equation}
Recall the Taylor's formula
\begin{equation}\label{Taylor}
J_j(\zeta)=J_j(\xi)+\sum_{k=1}^n(\zeta_k-\xi_k)\int_0^1\partial_{\xi_k}J_j(\xi+t(\zeta-\xi))\dd t,
\end{equation}
where
\begin{equation}
\frac{\partial}{\partial\xi_k}J_j(\xi)=|\xi|^{p-2}\left(\delta_{jk}+(p-2)\frac{\xi_j\xi_k}{|\xi|^2}\right).
\end{equation}
Then
\begin{equation}
\nabla\cdot(|\nabla v_\tau|^{p-2}\nabla v_\tau)=\nabla\cdot(|\nabla v_0|^{p-2}\nabla v_0)+\tau\nabla\cdot(\mathcal{A}^\tau\nabla V_\tau),
\end{equation}
where $\mathcal{A}^\tau$ is the matrix with coefficients
\begin{equation}
\mathcal{A}^\tau_{jk}=\int_0^1\pr_{\xi_k}J_j(\nabla v_0+t\tau\nabla V_\tau)\dd t.
\end{equation}
Therefore $V_\tau$ satisfies the equation
\begin{equation}
\left\{\begin{array}{l}\nabla\cdot(\mathcal{A}^\tau\nabla V_\tau)=0,\\[5pt] V_\tau|_{\pr\dom}=\phi_1.\end{array}\right.
\end{equation}
Since $v_0$ doesn't have critical points and $\tau V_\tau\to 0$ in $C^1(\dom)$, it follows that there is a $\tau_0>0$ such that the coefficient matrices $\mathcal{A}^\tau\in C^{0,\beta}(\dom)$ are elliptic, and have ellipticity bounds and $C^{0,\beta}$ bounds that are independent of $\tau\in(-\tau_0,\tau_0)$. By \cite[Theorems 8.34 \& 8.33]{GiTru}, $V_\tau$ is bounded in $C^{1,\beta}(\dom)$, uniformly for $\tau\in(-\tau_0,\tau_0)$.

By the theorem of Arzel\'a-Ascoli, there must exist a $V\in C^{1,\beta}(\dom)$ such that (after passing to a subsequence) $V_\tau\to V$ in $C^1(\dom)$. Clearly, this $V$ must be a weak solution to
\begin{equation}\label{eq-V}
\left\{\begin{array}{l}\nabla\cdot(A(v_0)\nabla V)=0,\\[5pt] V|_{\pr\dom}=\phi_1,\end{array}\right.
\end{equation}
where $A(v_0)$ is the matrix defined by \eqref{eq-def-A}.
Since the solutions to \eqref{eq-V} are unique, it follows that the $C^1$ convergence $V_\tau\to V$ holds even without passing to a subsequence. 
This shows that
\begin{equation}\label{eq-Frechet}
\left.C^1(\dom)-\frac{\dd}{\dd\tau}\right|_{\tau=0}v_\tau=V.
\end{equation}
Note here also the useful fact that if $\phi_1=\phi_0$ then $V=v_0$, so $v_0$ also solves
\begin{equation}
\nabla\cdot(A(v_0)\nabla v_0)=0.
\end{equation}

\section{Small data asymptotics when $1<p<2$}\label{small-asymptotics}

For the whole of this section, assume $p\in(1,2)$. As already discussed in the introduction, for small Dirichlet data it is reasonable to expect that the nonlinear term in equation \eqref{eq} dominates. This is made precise in this section, where small data asymptotics for the solutions of equation \eqref{eq} are derived. 

Suppose $v\in C^\infty(\dom)$ is a solution to
\begin{equation}\label{eq-p-Laplace1}
\nabla\cdot(|\nabla v|^{p-2}\nabla v)=0,
\end{equation}
without any critical points in $\dom$. As mentioned above, there are many such solutions.

For $1>\epsilon>0$, let $u_\epsilon\in C^{1,\beta}(\dom)$ be the solution of the boundary value problem
\begin{equation}\label{eq-epsilon}
\left\{\begin{array}{l}\nabla\cdot (\gamma(x)\nabla u_\epsilon)+\nabla\cdot(|\nabla u_\epsilon|^{p-2}\nabla u_\epsilon)=0,\\[5pt] u_\epsilon|_{\pr\dom}=\epsilon v|_{\pr\dom}.\end{array}\right.
\end{equation}
Let $R_\epsilon$ be defined via the Ansatz
\begin{equation}\label{Ansatz1}
u_\epsilon=\epsilon v+\epsilon^{3-p}R_\epsilon.
\end{equation}
By Theorem \ref{thm-existence},  $w_\epsilon=\epsilon^{2-p}R_\epsilon$ is bounded in $C^{1,\beta}(\dom)$. By the theorem of Arzel\'a-Ascoli, this implies that there exists $w_0\in C^{1,\beta}(\dom)$, such that (after passing to a subsequence) $w_\epsilon\to w_0$ in $C^{1}(\dom)$, as $\epsilon\to0$. Then
\begin{equation}
\mathscr{D}'(\dom)-\lim_{\epsilon\to0^+}\epsilon^{2-p}\left(\nabla\cdot (\gamma(x)\nabla u_\epsilon)+\nabla\cdot(|\nabla u_\epsilon|^{p-2}\nabla u_\epsilon)\right)=\nabla\cdot(|\nabla\tilde v|^{p-2}\nabla\tilde v),
\end{equation}
where $\tilde v=v+w_0$. It is then easy to see that $\tilde v$ is a weak solution to the equation
\begin{equation}
\left\{\begin{array}{l}\nabla\cdot(|\nabla \tilde v|^{p-2}\nabla \tilde v)=0,\\[5pt] \tilde v|_{\pr\dom}= v|_{\pr\dom}.\end{array}\right.
\end{equation}
Since solutions are unique, it follows that $w_0=0$. As the limit does not depend on the subsequence chosen,  then $\epsilon^{2-p}R_\epsilon\to 0$ in $C^{1}(\dom)$ without the need to pass to a subsequence.

By Taylor's formula,
\begin{equation}
\nabla\cdot(|\nabla u_\epsilon|^{p-2}\nabla u_\epsilon)=\epsilon^{p-1}\nabla\cdot(|\nabla v|^{p-2}\nabla v)
+\epsilon\nabla\cdot (A^\epsilon\nabla R_\epsilon),
\end{equation}
where
\begin{equation}
A^\epsilon_{jk}=\int_0^1\pr_{\xi_k}J_j\left(\nabla v+t\epsilon^{2-p}\nabla R_\epsilon\right)\dd t.
\end{equation}
Therefore, $R_\epsilon$ satisfies the equation
\begin{equation}\label{eq-A-epsilon}
\left\{\begin{array}{l} \nabla\cdot\left[(A^\epsilon+\epsilon^{2-p}\gamma)\nabla R_\epsilon\right]=-\nabla\cdot(\gamma\nabla v)\\[5pt] R_\epsilon|_{\pr\dom}=0.\end{array}\right.
\end{equation}
Since $\epsilon^{2-p}R_\epsilon\to 0$ in $C^{1}(\dom)$,  there exists $\epsilon_0>0$ so that the coefficient matrices $A^\epsilon+\epsilon^{2-p}\gamma\in C^{0,\beta}(\dom)$ are elliptic, and have ellipticity bounds and $C^{0,\beta}$ bounds that are independent of $\epsilon\in(0,\epsilon_0)$. By \cite[Theorems 8.34 \& 8.33]{GiTru}, $R_\epsilon$ is bounded in $C^{1,\beta}(\dom)$, uniformly for $\epsilon\in(0,\epsilon_0)$.

Applying again the theorem of Arze\'a-Ascoli, it follows that (after passing to a subsequence) $R_\epsilon\to R$, in $C^{1}(\dom)$.  Taking the limit $\epsilon\to0$ in \eqref{eq-A-epsilon} gives that  $R$ is a  solution to 
\begin{equation}\label{eq-R-1}
\left\{\begin{array}{l}\nabla\cdot(A(v)\nabla R)=-\nabla\cdot(\gamma\nabla v),\\[5pt] R|_{\pr\dom}=0,\end{array}\right.
\end{equation}
where $A(v)$ is the same matrix defined in \eqref{eq-def-A}.
Since the solution is unique, $R_\epsilon\to R$ in $C^1(\dom)$, without the need to pass to a subsequence.

Finally, a note on Dirichlet-to-Neumann maps. Denoting by $\Lambda_0$ the Dirichlet-to-Neumann map associated to equation \eqref{eq-p-Laplace}, the asymptotic expansion \eqref{Ansatz1} gives
\begin{equation}
\Lambda_\gamma(\epsilon v|_{\pr\dom})=\left.\left(\epsilon\nu\cdot\gamma\nabla v+\epsilon^{3-p}\nu\cdot\gamma\nabla R_\epsilon+\epsilon^{p-1}|\nabla v|\pr_\nu v+\epsilon\nu\cdot A^\epsilon\nabla R_\epsilon  \right)\right|_{\pr\dom},
\end{equation}
and therefore
\begin{equation}\label{eq-DN-expansion1}
C(\pr\dom)-\lim_{\epsilon\to0^+}\epsilon^{-1}\left(\Lambda_\gamma(\epsilon v)-\Lambda_0(\epsilon v)\right)
=\left.\nu\cdot\left(\gamma\nabla v+ A(v)\nabla R  \right)\right|_{\pr\dom}.
\end{equation}

\section{Large data asymptotics when $2<p<\infty$}\label{large-asymptotics} 

In this section the assumption is that $p\in(2,\infty)$. In this case the nonlinear term of equation \eqref{eq} should dominate at large Dirichlet data. Analogously to the previous section, large data asymptotics for the solutions of equation \eqref{eq} are derived below. As the arguments are very similar to those of section \ref{small-asymptotics}, they are presented in  a slightly abbreviated form.

Suppose $v\in C^{1,\beta}(\dom)$ is a solution to
\begin{equation}
\nabla\cdot(|\nabla v|^{p-2}\nabla v)=0,
\end{equation}
without critical points, and for $1>\epsilon>0$, let $u_\epsilon\in C^{1,\beta}(\dom)$ be the solution of the boundary value problem
\begin{equation}\label{eq-epsilon-2}
\left\{\begin{array}{l}\nabla\cdot (\gamma(x)\nabla u_\epsilon)+\nabla\cdot(|\nabla u_\epsilon|^{p-2}\nabla u_\epsilon)=0,\\[5pt] u_\epsilon|_{\pr\dom}=\epsilon^{-1}v|_{\pr\dom}.\end{array}\right.
\end{equation}
Let $\tilde R_\epsilon$ be defined via the Ansatz
\begin{equation}\label{Ansatz2}
u_\epsilon=\epsilon^{-1}v+\epsilon^{p-3} \tilde R_\epsilon.
\end{equation}
By essentially the same argument as in section \ref{small-asymptotics}, it follows that $\epsilon^{p-2} \tilde R_\epsilon\to 0$ in $C^{1}(\dom)$, as $\epsilon\to 0$. Also as above, $\tilde R_\epsilon$ satisfies the equation
\begin{equation}\label{eq-A-epsilon2}
\left\{\begin{array}{l} \nabla\cdot\left[(\tilde A^\epsilon+\epsilon^{p-2}\gamma)\nabla \tilde R_\epsilon\right]=-\nabla\cdot(\gamma\nabla v)\\[5pt] \tilde R_\epsilon|_{\pr\dom}=0.\end{array}\right.
\end{equation}
where
\begin{equation}
\tilde A^\epsilon_{jk}=\int_0^1\pr_{\xi_k} J_j(\nabla v+t\epsilon^{p-2}\nabla\tilde R_\epsilon)\dd t.
\end{equation}
There exists $\epsilon_0>0$ so that the coefficient matrices $\tilde A^\epsilon+\epsilon^{2-p}\gamma\in C^{0,\beta}(\dom)$ are elliptic, and have ellipticity bounds and $C^{0,\beta}$ bounds that are independent of $\epsilon\in(0,\epsilon_0)$. Therefore $\tilde R_\epsilon$ is bounded in $C^{1,\beta}(\dom)$, uniformly for $\epsilon\in(0,\epsilon_0)$, and in fact converges to a function $R\in C^{1,\beta}(\dom)$ which is the unique solution to  
\begin{equation}\label{eq-R-2}
\left\{\begin{array}{l}\nabla\cdot(A(v)\nabla R)=-\nabla\cdot(\gamma\nabla v),\\[5pt] R|_{\pr\dom}=0.\end{array}\right.
\end{equation}
Finally
\begin{equation}\label{eq-DN-expansion2}
C(\pr\dom)-\lim_{\epsilon\to0^+}\epsilon\left(\Lambda_\gamma(\epsilon^{-1} v)-\Lambda_0(\epsilon^{-1} v)\right)
=\left.\nu\cdot\left(\gamma\nabla v+ A(v)\nabla R  \right)\right|_{\pr\dom}.
\end{equation}

\section{Proof of Theorem \ref{thm-main}}\label{proof-main-thm}

For all values of $p$, and for each choice of $p$-harmonic function $v$ without critical points, the Dirichlet-to-Neumann map 
$\Lambda_\gamma$ determines the quantity $\left.\nu\cdot\left(\gamma\nabla v+ A(v)\nabla R  \right)\right|_{\pr\dom}$ via either equation \eqref{eq-DN-expansion1} or equation \eqref{eq-DN-expansion2}. Multiplying by any function $W\in C^\infty(\dom)$ and integrating over the boundary gives
\begin{multline}
I(v,W):=\int_{\pr\dom} \nu\cdot\left(\gamma\nabla v+ A(v)\nabla R  \right) W\dd x\\[5pt]
=\int_\dom \left[W\left(\nabla\cdot(\gamma\nabla v+A(v)\nabla R\right)+\nabla W\cdot\left(\gamma\nabla v+A(v)\nabla R\right)  \right]\dd x\\[5pt]
=\int_\dom \nabla v\cdot\gamma\nabla W\dd x-\int_\dom R\nabla\cdot(A(v)\nabla W)\dd  x.
\end{multline}

Let $\tau\in(-1,1)$ be a small parameter, and $v_0$ a smooth $p$-harmonic function without critical points. With the notation of section \ref{p-Laplace-linearization}, let $v=v_\tau$ and let $W$ be a solution to $\nabla\cdot(A(v_0)\nabla W)=0$. Then
\begin{multline}\label{eq-integral-identity}
\lim_{\tau\to0}\frac{1}{\tau}\left[I(v_\tau,W)-I(v_0,W)\right]=\int_\dom \nabla V\cdot\gamma\nabla W\dd x\\[5pt] -\int_\dom\dot R\nabla\cdot(A(v_0)\nabla W)\dd x
-\int_\dom R\nabla\cdot (\dot A(v_0,V)\nabla W)\dd x\\[5pt]
=\int_\dom \nabla V\cdot\gamma\nabla W\dd x-\int_\dom R\nabla\cdot (\dot A(v_0,V)\nabla W)\dd x,
\end{multline}
where $\dot R=\left.\frac{\dd}{\dd\tau}\right|_{\tau=0} R$, and
\begin{multline}
\dot A(v_0, V)=(p-2)|\nabla v_0|^{p-4}\Big[(\nabla v_0\cdot\nabla V)I\\[5pt]+(p-4)(\nabla v_0\cdot\nabla V)\frac{\nabla v_0\otimes\nabla v_0}{|\nabla v_0|^2} +\nabla v_0\otimes\nabla V+\nabla V\otimes\nabla v_0  \Big].
\end{multline}

Let $z\in\R^n$, $|z|^2=1$, and set $v_0=z\cdot x$. Then
\begin{equation}
A(v_0)=I+(p-2)z\otimes z
\end{equation}
and
\begin{multline}
\dot A(v_0, V)=(p-2)\Big[(z\cdot\nabla V)I\\[5pt]+(p-4)(z\cdot\nabla V) z\otimes z +z\otimes\nabla V+\nabla V\otimes z  \Big].
\end{multline}
The corresponding $R$ term satisfies
\begin{equation}\label{eq-R-z}
\left\{\begin{array}{l}\nabla\cdot\left[(I+(p-2)z\otimes z)\nabla R\right]=-\nabla\cdot(\gamma z),\\[5pt] R|_{\pr\dom}=0.\end{array}\right.
\end{equation}

\subsubsection*{The first set of special solutions}

Suppose $z,\xi,\eta\in\R^n\setminus\{0\}$, $t,s\in\R$, are such that $z\perp\xi,\eta$, $\xi\perp\eta$, $|z|=|\eta|=1$. Let 
\begin{equation}
\zeta_\pm=\pm s z+i(\xi\pm t\eta).
\end{equation}
 Then
\begin{equation}
\zeta_\pm\cdot A(v_0)\zeta_\pm=(p-1)s^2-\xi^2-t^2,
\end{equation}
therefore the choice 
\begin{equation}
s=(p-1)^{-\frac{1}{2}}\sqrt{\xi^2+t^2}
\end{equation}
 insures that both $e^{\zeta_\pm\cdot x}$ are solutions to equation $\nabla\cdot(A(v_0)\nabla V)=0$. Note that
\begin{equation}
\zeta_++\zeta_-=2i\xi.
\end{equation}
With the choices 
\begin{equation}
V=e^{\zeta_+\cdot x}, \quad W=e^{\zeta_-\cdot x},
\end{equation}
 it holds that
\begin{equation}\label{eq-V-W}
\int _\dom\nabla V\cdot\gamma\nabla W\dd x=-t^2\left(\frac{s}{t}z+i\eta\right)\cdot\hat\gamma(2\xi)\left(\frac{s}{t}z+i\eta\right)
-\xi\cdot\hat\gamma(2\xi)\xi.
\end{equation}
Note that 
\begin{equation}
\frac{s}{t}=(p-1)^{-\frac{1}{2}}\sqrt{1+\frac{\xi^2}{t^2}}
=(p-1)^{-\frac{1}{2}}\left( 1+\frac{\xi^2}{2}t^{-2} -\frac{|\xi|^4}{8}t^{-4}\right)+\Ord(t^{-6})
\end{equation}
as $t\to\infty$, and so
\begin{multline}
\int _\dom\nabla V\cdot\gamma\nabla W\dd x
=-t^2\left(\frac{z}{\sqrt{p-1}}+i\eta\right)\cdot\hat\gamma(2\xi)\left(\frac{z}{\sqrt{p-1}}+i\eta\right)\\[5pt]
-i\frac{\xi^2}{\sqrt{p-1}}\eta\cdot\hat\gamma(2\xi)z-\xi\cdot\hat\gamma(2\xi)\xi\\[5pt]
+it^{-2}\frac{|\xi|^4}{4\sqrt{p-1}}\eta\cdot\hat\gamma(2\xi)z 
+\Ord(t^{-4}).
\end{multline}
On the other hand
\begin{multline}
\int_\dom R\nabla\cdot (\dot A(v_0,V)\nabla W)\dd x=
2i(p-2)\hat R(2\xi)\Big[(z\cdot\zeta_+)(\xi\cdot\zeta_-)\\[5pt]
+(p-4)(z\cdot\zeta_+)(z\cdot\zeta_-)(z\cdot\xi)+(\zeta_+\cdot\zeta_-)(z\cdot\xi)+(z\cdot\zeta_-)(\xi\cdot\zeta_+)
\Big]\\
=0.
\end{multline}

Each order in $t$ in \eqref{eq-integral-identity}, as $t\to\infty$, can be computed separately from $\Lambda_\gamma$. That is, one can explicitly compute
\begin{equation}\label{eq-order-0}
\frac{z\cdot\hat\gamma(2\xi)z}{p-1}-\eta\cdot\hat\gamma(2\xi)\eta+\frac{2i}{\sqrt{p-1}}\eta\cdot\hat\gamma(2\xi)z,
\end{equation}
\begin{equation}\label{eq-order-1}
|\xi|^{-2}\xi\cdot\hat\gamma(2\xi)\xi+i(p-1)^{-\frac{1}{2}}\eta\cdot\hat\gamma(2\xi)z,
\end{equation}
\begin{equation}\label{eq-order-2}
\eta\cdot\hat\gamma(2\xi)z,
\end{equation}
and so on. Combining \eqref{eq-order-2} and \eqref{eq-order-0}, it follows that $\Lambda_\gamma$ determines 
\begin{equation}
\frac{z\cdot\hat\gamma(2\xi)z}{p-1}-\eta\cdot\hat\gamma(2\xi)\eta.
\end{equation}
Since the roles of $z$ and $\eta$ can be reversed, it follows that $\Lambda_\gamma$ also determines
\begin{equation}
\frac{\eta\cdot\hat\gamma(2\xi)\eta}{p-1}-z\cdot\hat\gamma(2\xi)z,
\end{equation}
and since $p\neq2$, it is possible to solve for
\begin{equation}
\eta\cdot\hat\gamma(2\xi)\eta.
\end{equation}
By combining \eqref{eq-order-2} and \eqref{eq-order-0}, it follows that $\Lambda_\gamma$ also determines
\begin{equation}
\xi\cdot\hat\gamma(2\xi)\xi.
\end{equation}

In summary, it has now been established that $\Lambda_\gamma$ determines
  the matrix elements
\begin{equation}
\xi\cdot\hat\gamma(2\xi)\xi,\quad \eta\cdot\hat\gamma(2\xi)\eta,\quad \mu\cdot\hat\gamma(2\xi)\eta, \quad \forall \xi,\mu,\eta\in\R^n\setminus\{0\}, \mu,\nu\perp\xi, \mu\perp\nu.
\end{equation}
In order to show that $\Lambda_\gamma$ determines $\gamma$, it is sufficient to further show that it determines matrix elements of the form 
\begin{equation}
\eta\cdot\hat\gamma(2\xi)\xi,\quad \forall \eta\perp\xi.
\end{equation}

\subsubsection*{The second set of special solutions}

Suppose $z,\xi,\eta,\mu\in\R^n\setminus\{0\}$, $t,s\in\R$, are such that $\eta\perp\xi,\mu,z$, $|z|=|\eta|=|\mu|=1$, $\mu$ and $\xi$ make an angle of $2\theta$, and $z\in\spn\{\mu,\xi\}$ bisects the angle between the two, i.e. $z$ makes an angle of $\theta$ with both $\xi$ and $\mu$.

Let 
\begin{equation}
\zeta_\pm=\pm s\mu+i(\xi\pm t\eta).
\end{equation}
Then
\begin{multline}
\zeta_\pm\cdot A(v_0)\zeta_\pm
=s^2-\xi^2-t^2\pm2is\mu\cdot\xi\\[5pt]
+(p-2)\left[s^2(\mu\cdot z)^2-(\xi\cdot z)^2\pm2is(\mu\cdot z)(\xi\cdot z)  \right]
\end{multline}
It follows that $e^{\zeta_\pm\cdot x}$ are solutions of $\nabla\cdot(A(v_0)\nabla V)=0$ if
\begin{equation}\label{eq-angle-condition-a}
\mu\cdot\xi+(p-2)(\mu\cdot z)(\xi\cdot z)=0
\end{equation}
and
\begin{equation}\label{eq-s-condition-a}
s^2\left[1+(p-2)(\mu\cdot z)^2  \right]=t^2+\xi^2\left[1+(p-2)\frac{(\xi\cdot z)^2}{|\xi|^2} \right].
\end{equation}
The identity \eqref{eq-angle-condition-a} is satisfied if $\theta\in\left(0,\frac{\pi}{2}\right)$ is chosen so that
\begin{equation}
\cos(2\theta)=\frac{2}{p}-1,
\end{equation}
in which case 
\begin{equation}
\cos^2\theta=\frac{1}{2}\left(1+\cos(2\theta)\right)=\frac{1}{p}.
\end{equation}
The identity \eqref{eq-s-condition-a} is satisfied by making the choice
\begin{equation}
s=\left[1+(p-2)(\mu\cdot z)^2  \right]^{-\frac{1}{2}}\left[t^2+\xi^2+(p-2)(\xi\cdot z)^2  \right]^{\frac{1}{2}}.
\end{equation}
As above, it is useful to note the expansion
\begin{multline}
\frac{s}{t}=\left[1+(p-2)(\mu\cdot z)^2  \right]^{-\frac{1}{2}}\\[5pt]
\times\left(1+\frac{1}{2}\left[ \xi^2+(p-2)(\xi\cdot z)^2 \right]t^{-2}
-\frac{1}{8}\left[ \xi^2+(p-2)(\xi\cdot z)^2 \right]^2t^{-4}  \right)\\[5pt]+\Ord(t^{-6}).
\end{multline}

With the choices
\begin{equation}
V=e^{\zeta_+\cdot x}, \quad W=e^{\zeta_-\cdot x},
\end{equation}
 it holds that
\begin{multline}
\int_\dom \nabla V\cdot\gamma\nabla W\dd x=
-t^2\left(\frac{s}{t}\mu+i\eta\right)\cdot\hat\gamma(2\xi)\left(\frac{s}{t}\mu+i\eta\right)
-\xi\cdot\hat\gamma(2\xi)\xi\\[5pt]
= -t^2\left(\frac{\mu}{\sqrt{1+(p-2)(\mu\cdot z)^2}}+i\eta\right)\cdot\hat\gamma(2\xi)\left(\frac{\mu}{\sqrt{1+(p-2)(\mu\cdot z)^2}}+i\eta\right)\\[5pt]
-i\frac{ \xi^2+(p-2)(\xi\cdot z)^2}{\sqrt{1+(p-2)(\mu\cdot z)^2}}\eta\cdot\hat\gamma(2\xi)\mu-\xi\cdot\hat\gamma(2\xi)\xi\\[5pt]
+i t^{-2}\frac{\left[ \xi^2+(p-2)(\xi\cdot z)^2\right]^2}{4\sqrt{1+(p-2)(\mu\cdot z)^2}}\eta\cdot\hat\gamma(2\xi)\mu
+\Ord(t^{-4}).
\end{multline}

On the other hand
\begin{multline}
\int_\dom R\nabla\cdot(\dot A(v_0,V)\nabla W)\dd x
=-2i(p-2)\hat R(2\xi)\Big[2s^2(z\cdot\mu)(\xi\cdot\mu)\\[5pt]+2\xi^2(z\cdot\xi) +(z\cdot\xi)(s^2+\xi^2-t^2)
+(p-4)(x\cdot\xi)(s^2(z\cdot\mu)^2+(z\cdot\xi)^2)  \Big]\\[5pt]
=-2i(p-2)\hat R(2\xi)(z\cdot\xi)\\[5pt]
\times\left[ 2s^2\left(\mu\cdot\frac{\xi}{|\xi|}\right)+(p-4)s^2(z\cdot\mu)^2-t^2+3\xi^2+(p-4)(z\cdot\xi)^2
\right].
\end{multline}
Note that, however complicated the above expression, it only contains terms with powers $t^2$ and $t^0$. Identifying the coefficient of $t^{-2}$ in \eqref{eq-integral-identity} gives that $\Lambda_\gamma$ determines
\begin{equation}
\eta\cdot\hat\gamma(2\xi)\mu.
\end{equation}
Since $p\neq2$, the angle $2\theta\neq\pi/2$ so $\mu\not\perp\xi$. Therefore, for any $\omega\in\R^n\setminus\{0\}$ such that  $\omega\perp\xi,\eta$, there exist $a,b\in\R\setminus\{0\}$ such that
\begin{equation}
\mu=a\xi+b\omega
\end{equation}
satisfies the conditions above. As
\begin{equation}
\eta\cdot\hat\gamma(2\xi)\mu=a\eta\cdot\hat\gamma(2\xi)\xi+b\eta\cdot\hat\gamma(2\xi)\omega,
\end{equation}
and as, by \eqref{eq-order-2}, $\eta\cdot\hat\gamma(2\xi)\omega$ is already known to be determined by $\Lambda_\gamma$, it follows that $\Lambda_\gamma$ determines 
\begin{equation}
\eta\cdot\hat\gamma(2\xi)\xi.
\end{equation}
This completes the proof of Theorem \ref{thm-main}.

\paragraph{Acknowledgments:}
C.C. was supported by NSTC grant number 112-2115-M-A49-002.

\bibliography{nonlinearity}
\bibliographystyle{plain}

\end{document}